\documentclass[12pt]{amsart}
\usepackage{amsmath, amsfonts, amssymb, amsthm,hyperref,mathtools}
\usepackage{bm}
\allowdisplaybreaks[4]
\def\({\bg(}
\def\){\bg)}

\def\Tr{{\rm Tr}}

\def\v{{\bm v}}
\def\u{{\bm u}}

\def\pmod #1{\ ({\rm{mod}}\ #1)}
\def\mod #1{\ {\rm mod}\ #1}
\def\Ack{\medskip\noindent {\bf Acknowledgments}}

\theoremstyle{plain}
\newtheorem{theorem}{Theorem}[section]
\newtheorem{lemma}{Lemma}

\theoremstyle{definition}

\theoremstyle{remark}
\newtheorem{remark}{Remark}

\makeatletter
\@namedef{subjclassname@2020}{%
	\textup{2020} Mathematics Subject Classification}
\makeatother
\vspace{4mm}

\begin{document}
	\medskip
	
	\title[Gaussian hypergeometric functions and cyclotomic matrices]
	{Gaussian hypergeometric functions and cyclotomic matrices}
	\author[H.-L. Wu and L.-Y. Wang]{Hai-Liang Wu and Li-Yuan Wang*}
	
	\address {(Hai-Liang Wu) School of Science, Nanjing University of Posts and Telecommunications, Nanjing 210023, People's Republic of China}
	\email{\tt whl@njupt.edu.cn}
	
    \address {(Li-Yuan Wang) School of Physical and Mathematical Sciences, Nanjing Tech University, Nanjing 211816, People's Republic of China}
    \email{\tt  lywang@njtech.edu.cn}
	
	\keywords{Jacobi sums, Gauss Sums, Finite Fields, Cyclotomic Matrices.
		\newline \indent This work was supported by the Natural Science Foundation of China (Grant Nos. 12101321 and 12201291).\newline \indent Corresponding author.}
	
	\begin{abstract}
		Let $q=p^n$ be an odd prime power and let $\mathbb{F}_q$ be the finite field with $q$ elements. Let $\widehat{\mathbb{F}_q^{\times}}$ be the group of all multiplicative characters of $\mathbb{F}_q$ and let $\chi$ be a generator of  $\widehat{\mathbb{F}_q^{\times}}$. In this paper, we investigate arithmetic properties of certain cyclotomic matrices involving nonzero squares over $\mathbb{F}_q$. For example, let $s_1,s_2,\cdots,s_{(q-1)/2}$ be all nonzero squares over $\mathbb{F}_q$. For any integer $1\le r\le q-2$, define the matrix 
		$$B_{q,2}(\chi^r):=\left[\chi^r(s_i+s_j)+\chi^r(s_i-s_j)\right]_{1\le i,j\le (q-1)/2}.$$
		We prove that if $q\equiv 3\pmod 4$, then  
			$$\det (B_{q,2}(\chi^r))=\prod_{0\le k\le (q-3)/2}J_q(\chi^r,\chi^{2k})=
		\begin{cases}
			(-1)^{\frac{q-3}{4}}{\bf i}^nG_q(\chi^r)^{\frac{q-1}{2}}/\sqrt{q} & \mbox{if}\ r\equiv 1\pmod 2,\\	
			G_q(\chi^r)^{\frac{q-1}{2}}/q                                     & \mbox{if}\ r\equiv 0\pmod 2,
		\end{cases}$$
	where $J_q(\chi^r,\chi^{2k})$ and $G_q(\chi^r)$ are the Jacobi sum and the Gauss sum over $\mathbb{F}_q$ respectively.
	\end{abstract}
	\maketitle
	
	\section{Introduction}
	\setcounter{lemma}{0}
	\setcounter{theorem}{0}
	\setcounter{equation}{0}
	\setcounter{conjecture}{0}
	\setcounter{remark}{0}
	\setcounter{corollary}{0}

    \subsection{Notation}
     Throughout this paper, for any $n\times n$ matrix $M$ the symbol $\det(M)$ denotes the determinant of $M$. Let $q$ be an odd prime power and let $\mathbb{F}_q$ be the finite field with $q$ elements. Let $\mathbb{F}_q^{\times}$ be the cyclic group of all nonzero elements of $\mathbb{F}_q$ and let $\widehat{\mathbb{F}_q^{\times}}$ be the group of all multiplicative characters of $\mathbb{F}_q$. For any $\psi\in\widehat{\mathbb{F}_q^{\times}}$, we extend $\psi$ to $\mathbb{F}_q$ by defining $\psi(0)=0$. We use the symbols $\varepsilon$ and $\phi$ to denote the trivial character and the unique quadratic character of $\mathbb{F}_q$ respectively, i.e., 
     $$\phi(x)=
       \begin{cases}
       	1  & \mbox{if}\ x\ \text{is a nonzero square},\\
       	0  & \mbox{if}\ x=0,\\
       	-1 & \mbox{otherwise}.
       \end{cases}$$
     
     Also, for any $A,B\in\widehat{\mathbb{F}_q^{\times}}$, the Jacobi sum $J_q(A,B)$ is defined by 
     $$J_q(A,B)=\sum_{x\in\mathbb{F}_q}A(x)B(1-x).$$
     Let $p$ be the characteristic of $\mathbb{F}_q$ with  $[\mathbb{F}_q:\mathbb{F}_p]=n$. Set $\zeta_p=e^{2\pi{\bf i}/p}$, where ${\bf i}$ is the primitive $4$-th root of unity with argument $\pi/2$. Then the Gauss sum $G_q(A)$ over $\mathbb{F}_q$ is defined by
     $$G_q(A)=\sum_{x\in\mathbb{F}_q}A(x)\zeta_p^{\Tr(x)},$$
     where 
     $$\Tr(x):=\Tr_{\mathbb{F}_q/\mathbb{F}_p}(x)=\sum_{j=0}^{n-1}x^{p^j}$$
     is the trace map from $\mathbb{F}_q$ to $\mathbb{F}_p$. 
     
     \subsection{Background and Motivation} Let $p$ be an odd prime and let $(\frac{\cdot}{p})$ be the Legendre symbol. Then it is known that 
     $$\sum_{x\in\mathbb{F}_p}\left(\frac{x}{p}\right)\zeta_p^x=
     \begin{cases}
     	\sqrt{p}        & \mbox{if}\ p\equiv 1\pmod 4,\\
     	{\bf i}\sqrt{p} & \mbox{if}\ p\equiv 3\pmod 4.
     \end{cases}$$
     By this and the Hasse-Davenport lifting formula (see \cite[Theorem 3.7.4]{Cohen}), one can easily obtain the explicit value of $G_q(\phi)$, where $q=p^n$ is an odd prime power. In fact, we have (see \cite[Corollary 3.7.6]{Cohen})
     \begin{equation}\label{Eq. explicit value of quadratic gauss sum}
     	G_q(\phi)=\begin{cases}
     		(-1)^{n-1}\sqrt{q}          & \mbox{if}\ p\equiv 1\pmod 4,\\
     		(-1)^{n-1}{\bf i}^n\sqrt{q} & \mbox{if}\ p\equiv 3\pmod 4.
     	\end{cases}
     \end{equation}
 
     On the other hand, the arithmetic properties of cyclotomic matrices have been extensively investigated and these cyclotomic matrices have close relations with Gauss sums. For example, let $\psi$ be a nontrivial multiplicative character of $\mathbb{F}_p$. Carlitz \cite[Theorem 4]{Carlitz} studied the matrix 
     \begin{equation*}
     	C_p:=\left[\psi(j-i)\right]_{1\le i,j\le p-1}.
     \end{equation*}
     Surprisingly, the eigenvalues of $C_p$ are closely related to the Gauss sums over $\mathbb{F}_p$. In fact,   
     Carlitz proved that the characteristic polynomial of $C_p$ is given by 
     \begin{equation*}
     	f(x)=\det\left(xI_{p-1}-C_p\right)
     	=\left(x^f-G_p(\psi)^f\right)^{t-1}\left(x^f-\frac{1}{p}G_p(\psi)^f\right),
     \end{equation*}
     where $f=\min\{k\in\mathbb{Z}^+: \psi^k=\varepsilon\}$ is the order of $\psi$, and  $t=(p-1)/f$. Hence we have  
     \begin{equation}\label{Eq. det of Carlitz matrix}
     	\det (C_p)=(-1)^t\frac{G_p(\psi)^{p-1}}{p}.
     \end{equation}

     Along this line, Chapman \cite{Chapman} further investigated some variants of Carlitz's matrix $C_p$. For instance, let $\epsilon_p>1$ and $h_p$ be the fundamental unit and the class number of $\mathbb{Q}(\sqrt{p})$ respectively. Then Chapman's well-known ``evil determinant" conjecture states that 
     \begin{equation*}
     	\det\left[\left(\frac{j-i}{p}\right)\right]_{1\le i,j\le \frac{p+1}{2}}=
     	\begin{cases}
     		-a_p' & \mbox{if}\ p\equiv 1\pmod 4,\\
     		1     & \mbox{if}\ p\equiv 3\pmod 4,
     	\end{cases}
     \end{equation*}
     where $a_p'$ is defined by 
     $$\epsilon_p^{(2-(\frac{2}{p}))h_p}=a_p'+b_p'\sqrt{p}\ (a_p',b_p'\in\mathbb{Q}).$$ 
     This challenging conjecture was later confirmed by Vsemirnov \cite{Vsemirnov12,Vsemirnov13}. 
	
	In 2019, Z.-W. Sun initiated the study of the cyclotomic matrices involving quadratic polynomials over finite fields. For example, Sun \cite[Theorem 1.2]{Sun19} considered the matrix 
	$$S_p:=\left[\left(\frac{i^2+j^2}{p}\right)\right]_{1\le i,j\le (p-1)/2}.$$
	As a matrix over $\mathbb{Z}$, Sun conjectured that $-\det(S_p)$ is a square of some integer if  $p\equiv 3\pmod 4$. In the case $p\equiv 1\pmod4$, it is known that $p$ can be uniquely written as $p=a^2+4b^2$ with $a,b\in\mathbb{Z}$, $b>0$ and $a\equiv 1\pmod4$.  Cohen, Sun and Vsemirnov (see \cite[Remark 4.2]{Sun19}) conjectured that $\det (S_p)/a$ is also a square of some integer. Later, the case $p\equiv 3\pmod 4$ was confirmed by Alekseyev and Krachun, and the case $p\equiv 1\pmod 4$ was proved by the first author \cite{Wu-CR}. Moreover, the authors and Y.-F. She \cite{Wu-She-Wang} studied some generalizations of $S_p$. 
	
	Recently, Z.-W. Sun \cite{Sun24} posed many conjectures related to cyclotomic matrices. For example, Sun \cite[Theorem 1.1]{Sun24} studied the matrix 
	\begin{equation*}
		A_p(x)=\left[x+\left(\frac{i^2+j^2}{p}\right)+\left(\frac{i^2-j^2}{p}\right)\right]_{1\le i,j\le (p-1)/2}.
	\end{equation*}
	For $p\equiv 1\pmod 4$, Sun proved that $A_p(x)$ is a singular matrix. For $p\equiv 3\pmod4$, Sun \cite[Remark 1.1]{Sun24} conjectured that 
	$$\det (A_p(x))=\left(\frac{p-1}{2}x-1\right)p^{\frac{p-3}{4}}.$$
	This conjecture was later confirmed by Li and the first author \cite{Li-Wu}. 
	
	After introducing the above relevant research results, we now describe our research motivations. Gauss sums are closely related to Jacobi sums, and by (\ref{Eq. det of Carlitz matrix}) the determinant of Carlitz's matrix $C_p$ involves Gauss sums closely. Motivated by this, it is natural to search for a cyclotomic matrix whose determinant  is closely related to Jacobi sums. Fortunately, we will see below that a generalization of Sun's matrix $A_p(0)$ can meet our requirement. 
	
	Throughout the remaining part of this paper, we always let $\chi$ be a generator of $\widehat{\mathbb{F}_q^{\times}}$ and let $1\le r\le q-2$ be an integer. For any integer $1<k<q-1$ with $k\mid q-1$, let 
	$$D_k:=\{a_1,a_2,\cdots,a_{(q-1)/k}\}=\{x^k: x\in\mathbb{F}_q^{\times}\}$$
	be the set of all nonzero $k$-th powers over $\mathbb{F}_q$. As a generalization of Sun's matrix $A_p(0)$, we define the matrix 
	\begin{equation}\label{Eq. Wu-Wang matrix }
		B_{q,k}(\chi^r):=\left[\chi^r(a_i+a_j)+\chi^r(a_i-a_j)\right]_{1\le i,j\le (q-1)/k}.
	\end{equation}
	
	\subsection{The Gaussian hypergeometric function} To state our main results on $B_{q,k}(\chi^r)$, we first introduce some basic notations on Gaussian hypergeometric functions. 
	
	Greene \cite{Greene} initiated the study of Gaussian hypergeometric functions. Greene first used Jacobi sums to obtain finite field analogues of binomial coefficients. For any $A,B\in\widehat{\mathbb{F}_q^{\times}}$, Greene \cite[Definition 2.4]{Greene} defined $\binom{A}{B}$ by 
	$$\binom{A}{B}:=\frac{B(-1)}{q}J_q(A,\bar{B}).$$
	Readers may refer to \cite[(2.6)--(2.16)]{Greene} for the basic properties of $\binom{A}{B}$.	
	
	Now for characters $A_0,A_1,\cdots,A_n,B_1,\cdots,B_n\in\widehat{\mathbb{F}_q^{\times}}$, the Gaussian hypergeometric function from $\mathbb{F}_q$ to $\mathbb{C}$ (see \cite[Definition 3.10]{Greene}) is defined by 
	\begin{equation*}
		_{n+1}F_n\left(\begin{array}{cccccc}
		A_0 & A_1 & \cdots & A_n\\
		    & B_1 & \cdots & B_n
		\end{array}\Bigg|x\right)_q:=\frac{q}{q-1}\sum_{\chi\in\widehat{\mathbb{F}_q^{\times}}}\binom{A_0\chi}{\chi}\binom{A_1\chi}{B_1\chi}\cdots\binom{A_n\chi}{B_n\chi}\chi(x).
	\end{equation*}
    This can be viewed as a finite field analogue of the classical hypergeometric function
    \begin{equation*}
    	_{n+1}F_n\left(\begin{array}{cccccc}
    		a_0 & a_1 & \cdots & a_n\\
    		    & b_1 & \cdots & b_n
    	\end{array}\Bigg|x\right):=\lambda\sum_{k=0}^{+\infty}\binom{a_0+k-1}{k}\binom{a_1+k-1}{b_1+k-1}\cdots\binom{a_n+k-1}{b_n+k-1}x^k,
    \end{equation*}
     where
     $$\lambda=\prod_{k=1}^n\binom{a_k-1}{b_k-1}^{-1}.$$
     
	Gaussian hypergeometric functions are closely related to the numbers of rational points on varieties over finite fields and modular forms. Readers may refer to the survey paper \cite{F}.

	\subsection{Main results} We now state our main results. For $k=2$ we have the following theorem.
	
	\begin{theorem}\label{Thm. Bq(r)}
		Let $q=p^n$ be an odd prime power and let $\chi$ be a generator of $\widehat{\mathbb{F}_q^{\times}}$. Let $1\le r\le q-2$ be an integer. Then the following results hold.
		
		{\rm (i)} Suppose $q\equiv 3\pmod4$. Then 
		\begin{equation}\label{Eq. in Thm. 1}
			\det (B_{q,2}(\chi^r))=\prod_{0\le k\le (q-3)/2}J_q(\chi^r,\chi^{2k})=
			\begin{cases}
				(-1)^{\frac{q-3}{4}}{\bf i}^nG_q(\chi^r)^{\frac{q-1}{2}}/\sqrt{q} & \mbox{if}\ r\equiv 1\pmod 2,\\	
				G_q(\chi^r)^{\frac{q-1}{2}}/q                                     & \mbox{if}\ r\equiv 0\pmod 2.
			\end{cases}
		\end{equation}
		
		{\rm (ii)} For any divisor $1<k<q-1$ of $q-1$, if $(q-1)/k\equiv 0\pmod 2$, then $B_{q,k}(\chi^r)$ is a singular matrix. 
	\end{theorem}
	
	\begin{remark}
		By Theorem \ref{Thm. Bq(r)}(i) and Greene's finite field analogues of binomial coefficients,
		\begin{equation}\label{Eq. in Remark 1}
				\det (B_{q,2}(\chi^r))=q^{\frac{q-1}{2}}\prod_{k=0}^{(q-3)/2}\binom{\chi^r}{\chi^{2k}}
		\end{equation}
		for $q\equiv 3\pmod 4$. 
	\end{remark}
	
    For $k=4$, by Theorem \ref{Thm. Bq(r)}(ii) the matrix $B_{q,4}(\chi^r)$ is singular if $q\equiv 1\pmod 8$. For the case $q\equiv 5\pmod 8$ we have the following result. 
	\begin{theorem}\label{Thm. 4th powers}
		 Let $q=p^n$ be a prime power with $q\equiv 5\pmod 8$ and let $\chi$ be a generator of $\widehat{\mathbb{F}_q^{\times}}$. Let $1\le r\le q-2$ be an integer. Then
		 	\begin{equation}\label{Eq. in Thm. 2}
		 		\det(B_{q,4}(\chi^r))=(-1)^r\left(\frac{q}{2}\right)^{\frac{q-1}{4}}
		 		\prod_{k=0}^{(q-5)/4}\ _{2}F_{1}
		 		\left(\begin{array}{cc}
		 			\chi^{-r}  &  \chi^{4k}\\
		 			           &  \chi^{4k+r}
		 		\end{array}\Bigg|-1\right)_q.
	 		\end{equation}
	\end{theorem}
	
	\begin{remark}\label{Remark of Thm. 2}
		 (i) Similar to (\ref{Eq. in Remark 1}), by \cite[(2.8) and Theorem 4.16]{Greene} for $q\equiv 5\pmod 8$ we have 
		$$\det(B_{q,4}(\chi^r))=\left(\frac{q}{2}\right)^{\frac{q-1}{4}}
		\prod_{k=0}^{(q-5)/4}\left(\binom{\chi^r}{\chi^{2k}}+\binom{\chi^r}{\phi\chi^{2k}}\right).$$
		
		(ii) When $r=(q-1)/2$, if we set $q-1=4m$, then by Theorem \ref{Thm. 4th powers} we have 
		\begin{equation}\label{Eq. in remark 2}
		\det(B_{q,4}(\phi))=(-1)^r\left(\frac{q}{2}\right)^{\frac{q-1}{4}}
		\prod_{k=0}^{m-1}\ _{2}F_{1}
		\left(\begin{array}{cc}
			\phi  &  \chi^{4k}\\
			&  \phi\chi^{4k}
		\end{array}\Bigg|-1\right)_q.
		\end{equation}
		By Kummer's theorem (see \cite[p.9]{Bailey}) for $0\le k\le m-1$ we have 
		\begin{equation*}
			_{2}F_{1}
			\left(\begin{array}{cc}
				1/2  &  k/m\\
				     &  1-1/2+k/m
			\end{array}\Bigg|-1\right)
		=\frac{\Gamma\left(\frac{1}{2}+\frac{k}{m}\right)\Gamma\left(1+\frac{k}{2m}\right)}
		{\Gamma\left(1+\frac{k}{m}\right)\Gamma\left(\frac{1}{2}+\frac{k}{2m}\right)},
		\end{equation*}
		where $\Gamma(\cdot)$ is the Gamma function. Using the Gauss multiplication formula (see \cite[Theorem 1.5.2]{Andrews})
		$$\Gamma(mz)(2\pi)^{(m-1)/2}=m^{mz-\frac{1}{2}}\Gamma(z)\Gamma\left(z+\frac{1}{m}\right)\cdots\Gamma\left(z+\frac{m-1}{m}\right),$$
		we obtain 
		$$\prod_{k=0}^{m-1}\ _{2}F_{1}
		\left(\begin{array}{cc}
			1/2  &  k/m\\
			&  1-1/2+k/m
		\end{array}\Bigg|-1\right)=\frac{\Gamma(m/2)m^{m/2}}{(m-1)!}
	\cdot\prod_{k=0}^{m-1}\frac{\Gamma\left(1+\frac{k}{2m}\right)}
	{\Gamma\left(\frac{1}{2}+\frac{k}{2m}\right)}.$$
	Thus, in the case $r=(q-1)/2$, equality (\ref{Eq. in Thm. 2}) can be viewed as a finite field analogue of the above equality.
		
		(iii) By Theorem \ref{Thm. Bq(r)}(ii) we have $\det(B_{q,2k+1}(\chi^r))=0$ for any $0\le k\le (q-3)/2$. For $3\le k\le (q-3)/2$, getting a simple expression of $\det(B_{q,2k}(\chi^r))$ like the equations (\ref{Eq. in Thm. 1}) and (\ref{Eq. in Thm. 2}) seems difficult.
	\end{remark}
	
	As an application of Theorem \ref{Thm. Bq(r)} and Theorem \ref{Thm. 4th powers}, we next consider the matrix over $\mathbb{F}_q$ defined by 
	\begin{equation}\label{Eq. definition of T_q(r)}
		T_{q,k}(r):=\left[(a_i+a_j)^{q-1-r}+(a_i-a_j)^{q-1-r}\right]_{1\le i,j\le (q-1)/k},
	\end{equation}
	where $1\le r\le (q-2)$ is an integer and $1<k<q-1$ is a divisor of $q-1$. We have the following result.
	
	\begin{theorem}\label{Thm. Tq(r)}
	 Let $q=p^n$ be an odd prime power and let $1\le r\le q-2$ be an integer. Then the following results hold.
	 
	 {\rm (i)} Let $1<k<q-1$ be a divisor of $q-1$ with $(q-1)/k\equiv 0\pmod 2$. Then $T_{q,k}(r)$ is a singular matrix over $\mathbb{F}_q$. 
	 
	 {\rm (ii)} $T_{q,2}(r)$ is nonsingular if and only if $p\equiv 3\pmod 4, n=1$ and $r\in\{1,2\}$. Moreover, if $p\equiv 3\pmod4$ and $p>3$, then 
	 $$\det (T_{p,2}(1))=(-1)^{\frac{p+3+2h(-p)}{4}}\in\mathbb{F}_p,$$
	 where $h(-p)$ is the class number of the imaginary quadratic field $\mathbb{Q}(\sqrt{-p})$. Also, we have 
	 $$\det (T_{p,2}(2))=(-1)^{\frac{p+1}{4}}\in\mathbb{F}_p.$$
	 
	 {\rm (iii)} If $q\equiv 5\pmod 8$, then 
	 $$ \det(T_{q,4}(r))=
	 \left(\frac{1}{2}\right)^{\frac{q-1}{4}}\prod_{k=0}^{(q-5)/4}
	 \left(\binom{2k+r}{r}+\binom{2k+r+(q-1)/2}{r}\right)\in\mathbb{F}_p.$$
	\end{theorem}
	
	\subsection{Outline of the paper} We will prove Theorems \ref{Thm. Bq(r)}--\ref{Thm. 4th powers} in Section 2, and the proof of Theorem \ref*{Thm. Tq(r)} will be given in Section 3.

	\section{Proofs of Theorems \ref*{Thm. Bq(r)}--\ref*{Thm. 4th powers}}
		\setcounter{lemma}{0}
	\setcounter{theorem}{0}
	\setcounter{equation}{0}
	\setcounter{conjecture}{0}
	\setcounter{remark}{0}
	\setcounter{corollary}{0}

	We begin with the Hasse-Davenport product formula (see \cite[Theorem 3.7.3]{Cohen}). 
	
	\begin{lemma}\label{Lem. the HD product formula}
		Let $\rho\in\widehat{\mathbb{F}_q^{\times}}$ be a character of order $m>1$ and let $\psi\in\widehat{\mathbb{F}_q^{\times}}$. Then 
		\begin{equation*}
			\prod_{0\le k\le m-1}G_q(\psi\rho^k)=-\psi^{-m}(m)G_q(\psi^m)\prod_{0\le k\le m-1}G_q(\rho^k).
		\end{equation*}
	\end{lemma}
	
	The following results are some basic properties of Gauss sums and Jacobi sums. Readers may refer to \cite{BEK}.
 	
 	\begin{lemma}\label{Lem. basic properties of Gauss sums}
 		Let $\psi,\psi_1,\psi_2\in\widehat{\mathbb{F}_q^{\times}}$ be nontrivial characters. Then the following results hold. 
 		
 		{\rm (i)} If $\psi_1\psi_2\neq\varepsilon$, then 
 		$$J_q(\psi_1,\psi_2)=\frac{G_q(\psi_1)G_q(\psi_2)}{G_q(\psi_1\psi_2)}.$$
 		
 		{\rm (ii)} $G_q(\psi)G_q(\psi^{-1})=q\psi(-1)$. 
 		
 		{\rm (iii)} $J_q(\psi,\psi^{-1})=-\psi(-1)$. 
 	\end{lemma}

	We also need the following known result in linear algebra. 
	
	\begin{lemma}\label{Lem. an easy lemma in linear algebra}
		Let $m$ be a positive integer and let $M$ be an $m\times m$ complex matrix. Let $\lambda_1,\cdots,\lambda_m\in\mathbb{C}$, and let $\v_1,\cdots,\v_m\in\mathbb{C}^m$ be column vectors. Suppose that 
		$$M\v_i=\lambda_i\v_i$$
		for any $1\le i\le m$ and $\v_1,\cdots,\v_m$ are linearly independent over $\mathbb{C}$. Then $\lambda_1,\cdots,\lambda_m$ are exactly all the eigenvalues of $M$ (counting multiplicity). 
	\end{lemma}
	
	Now we are in a position to prove our first theorem.
	
	{\noindent{\bf Proof of Theorem \ref{Thm. Bq(r)}}.} Throughout this proof, let 
	$$D_2=\{x^2:x\in\mathbb{F}_q^{\times}\}=\{s_1,s_2,\cdots,s_{(q-1)/2}\}.$$
	In $\mathbb{F}_q[x]$ one can easily verify that the polynomial $x^{\frac{q-1}{2}}-1$ has the factorization
	\begin{equation*}
		x^{\frac{q-1}{2}}-1=\prod_{1\le j\le (q-1)/2}\left(x-s_j\right). 
	\end{equation*}
	This implies 
	\begin{equation}\label{Eq. product of all squares over finite fields}
		\prod_{1\le j\le (q-1)/2}s_j=(-1)^{\frac{q+1}{2}}.
	\end{equation}
	As $\chi$ is a generator of $\widehat{\mathbb{F}_q^{\times}}$, we have $\chi^r(-1)=(-1)^r$. By this and (\ref{Eq. product of all squares over finite fields}) we obtain 
	\begin{equation}\label{Eq. 1st step in the proof of Thm. 1}
		\det(B_{q,k}(\chi^r))=\prod_{1\le j\le (q-1)/2}\chi^r(s_j)\cdot \det( M_q(\chi^r))=
		(-1)^{\frac{(q+1)r}{2}}\det(M_q(\chi^r)),
	\end{equation}
	where 
	$$M_q(\chi^r):=\left[\chi^r\left(1+\frac{s_j}{s_i}\right)+\chi^r\left(1-\frac{s_j}{s_i}\right)\right]_{1\le i,j\le (q-1)/2}.$$
	
	(i) Suppose now $q\equiv 3\pmod 4$. Then for any integer $1\le k\le (q-1)/2$, set the column vector 
	$$\v_k=\left(\chi^k(s_1),\chi^k(s_2),\cdots,\chi^k(s_{(q-1)/2})\right)^T,$$
	and let the complex number 
	\begin{equation}\label{Eq. eigenvalues of Mq(r)}
		\lambda_k=\sum_{j=1}^{(q-1)/2}\chi^r(1+s_j)\chi^k(s_j)+\sum_{j=1}^{(q-1)/2}\chi^r(1-s_j)\chi^k(s_j).
	\end{equation}
	Then for any $1\le i\le (q-1)/2$, one can verify that 
	\begin{align*}
		&\sum_{j=1}^{(q-1)/2}\left(\chi^r(1+s_j/s_i)+\chi^r(1-s_j/s_i)\right)\chi^k(s_j)\\
	   =&\sum_{j=1}^{(q-1)/2}\left(\chi^r(1+s_j/s_i)+\chi^r(1-s_j/s_i)\right)\chi^k(s_j/s_i)\chi^k(s_i)\\
	   =&\sum_{j=1}^{(q-1)/2}\left(\chi^r(1+s_j)+\chi^r(1-s_j)\right)\chi^k(s_j)\chi^k(s_i)\\
	   =&\lambda_k\chi^k(s_i).
	\end{align*}
	This illustrates that for any $1\le k\le (q-1)/2$ we have 
	\begin{equation}\label{Eq. step 2 in the proof of Thm. 1}
		M_q(\chi^r)\v_k=\lambda_k\v_k.
	\end{equation}
	Note also that 
	$$\det\left[\chi^i(s_j)\right]_{1\le i,j\le (q-1)/2}=\pm \prod_{1\le i<j\le (q-1)/2}\left(\chi(s_j)-\chi(s_i)\right)\neq 0.$$
	Hence $\v_1,\cdots,\v_{(q-1)/2}$ are linearly independent over $\mathbb{C}$. This, together with (\ref{Eq. eigenvalues of Mq(r)}) and Lemma \ref{Lem. an easy lemma in linear algebra}, implies that $\lambda_1,\cdots,\lambda_{(q-1)/2}$ are exactly all the eigenvalues of $M_q(\chi^r)$. We next evaluate $\lambda_k$ for $1\le k\le (q-1)/2$. 
	
	{\bf Case I:} $k$ is even. 
	
	In this case, by (\ref{Eq. eigenvalues of Mq(r)}) we have 
	\begin{align*}
		\lambda_k
		&=\sum_{j=1}^{(q-1)/2}\chi^r(1+s_j)\chi^k(s_j)+\sum_{j=1}^{(q-1)/2}\chi^r(1-s_j)\chi^k(-s_j)\\
		&=\sum_{x\in\mathbb{F}_q}\chi^r(1+x)\chi^k(x)\\
		&=\sum_{x\in\mathbb{F}_q}\chi^r(1+x)\chi^k(-x).
	\end{align*}
	This implies that for any even integer $1\le k\le (q-1)/2$ we have 
	\begin{equation}\label{Eq. lambda k for k even}
		\lambda_k=J_q(\chi^r,\chi^k).
	\end{equation}

	{\bf Case II:} $k$ is odd.
	
	In this case, we first observe that 
	\begin{align*}
		 &2\sum_{j=1}^{(q-1)/2}\chi^r(1-s_j)\chi^k(-s_j)\\
		=&\sum_{x\in\mathbb{F}_q}\left(\varepsilon(x)+\phi(x)\right)\chi^r(1-x)\chi^k(-x)\\
		=&-\sum_{x\in\mathbb{F}_q}\chi^r(1-x)\chi^k(x)-\sum_{x\in\mathbb{F}_q}\chi^r(1-x)\chi^{k+\frac{q-1}{2}}(x)\\
		=&-\left(J_q(\chi^r,\chi^k)+J_q(\chi^r,\chi^{k+\frac{q-1}{2}})\right).
	\end{align*}
	By this one can verify that 
	\begin{align*}
		\lambda_k
		&=\sum_{j=1}^{(q-1)/2}\chi^r(1+s_j)\chi^k(s_j)+\sum_{j=1}^{(q-1)/2}\chi^r(1-s_j)\chi^k(-s_j)-2\sum_{j=1}^{(q-1)/2}\chi^r(1-s_j)\chi^k(-s_j)\\
		&=\sum_{x\in\mathbb{F}_q}\chi^r(1+x)\chi^k(x)+\left(J(\chi^r,\chi^k)+J(\chi^r,\chi^{k+\frac{q-1}{2}})\right)\\
		&=-J_q(\chi^r,\chi^k)+\left(J_q(\chi^r,\chi^k)+J_q(\chi^r,\chi^{k+\frac{q-1}{2}})\right)\\
		&=J_q(\chi^r,\chi^{k+\frac{q-1}{2}}).
	\end{align*}
	Hence for any odd integer $1\le k\le (q-1)/2$ we have 
	\begin{equation}\label{Eq. lambda k for k odd}
		\lambda_k=J(\chi^r,\chi^{k+\frac{q-1}{2}}).
	\end{equation} 
	Combining (\ref{Eq. 1st step in the proof of Thm. 1}) with (\ref{Eq. lambda k for k even}) and (\ref{Eq. lambda k for k odd}), we obtain 
	\begin{align}\label{Eq. product of Jacobi sums}
		\det(B_{q,k}(\chi^r))
		&=\prod_{1\le k\le (q-3)/4}J_q(\chi^r,\chi^{2k})\cdot\prod_{1\le k\le (q+1)/4}J_q(\chi^r,\chi^{2k-1+\frac{q-1}{2}}) \notag \\
		&=\prod_{0\le k\le (q-3)/2}J_q(\chi^r,\chi^{2k}).
	\end{align}
	
	We next compute the explicit value of (\ref{Eq. product of Jacobi sums}). We first consider the case $r\equiv1\pmod2$. In this case, $\chi^r\chi^{2k}\neq\varepsilon$ for any $0\le k\le (q-3)/2$. Hence by Lemma \ref{Lem. the HD product formula} and Lemma \ref{Lem. basic properties of Gauss sums} one can verify that 
	\begin{align*}
		  \prod_{0\le k\le (q-3)/2}J_q(\chi^r,\chi^{2k})
		&=\prod_{0\le k\le (q-3)/2}\frac{G_q(\chi^r)G_q(\chi^{2k})}{G_q(\chi^{r+2k})}\\
		&=G_q(\chi^r)^{\frac{q-1}{2}}\prod_{0\le k\le (q-3)/2}\frac{G_q(\chi^{2k})}{G_q(\chi^{r+2k})}\\
		&=-G_q(\chi^r)^{\frac{q-1}{2}}\chi^{\frac{r(q-1)}{2}}\left(\frac{q-1}{2}\right)/G_q(\chi^{\frac{r(q-1)}{2}}).
	\end{align*}
	As $r$ is odd, we have $\chi^{\frac{r(q-1)}{2}}=\phi$. Also, note that 
	$$-\phi\left(\frac{q-1}{2}\right)=\phi(2)=(-1)^{\frac{q+1}{4}}.$$
	By the above and (\ref{Eq. explicit value of quadratic gauss sum}) we obtain that 
	\begin{equation*}
		\prod_{0\le k\le (q-3)/2}J_q(\chi^r,\chi^{2k})=(-1)^{\frac{q-3}{4}}{\bf i}^nG_q(\chi^r)^{\frac{q-1}{2}}/\sqrt{q}
	\end{equation*}
	for any odd integer $1\le r\le (q-2)$. 
	
	For the case $r\equiv 0\pmod 2$, by Lemma \ref{Lem. the HD product formula} and Lemma \ref{Lem. basic properties of Gauss sums} one can verify that 
	\begin{align*}
		 \prod_{0\le k\le (q-3)/2}J_q(\chi^r,\chi^{2k})
		&=J_q(\chi^r,\chi^{-r})\cdot\frac{G_q(\varepsilon)}{G_q(\chi^r)G_q(\chi^{-r})}\cdot\prod_{0\le k\le (q-3)/2}\frac{G_q(\chi^r)G_q(\chi^{2k})}{G_q(\chi^{r+2k})}\\
		&=-\frac{1}{q}G_q(\chi^r)^{\frac{q-1}{2}}\chi^{\frac{r(q-1)}{2}}\left(\frac{q-1}{2}\right)/G_q(\chi^{\frac{r(q-1)}{2}}).
	\end{align*}
	As $r\equiv0\pmod2$, the character $\chi^{\frac{r(q-1)}{2}}=\varepsilon$ and hence 
	\begin{equation*}
	\prod_{0\le k\le (q-3)/2}J_q(\chi^r,\chi^{2k})=G_q(\chi^r)^{\frac{q-1}{2}}/q.     
	\end{equation*} 
	This completes the proof of (i). 
	
	(ii) As $q\equiv 1\pmod{2k}$, we have $-1\in D_k$. Hence for any $1\le j\le (q-1)/k$, there is a $j'$ such that $1\le j'\le (q-1)/2, j'\neq j$ and $a_j=-a_{j'}$. This implies that in the matrix $B_{q,k}(\chi^r)$ the $j$-th column is the same as the $j'$-th column. By this, $B_{q,k}(\chi^r)$ is singular if $(q-1)/k\equiv 0\pmod 2$. 
	
	In view of the above, we have completed the proof of Theorem \ref{Thm. Bq(r)}.\qed 
	
	We next prove our second theorem.
	
	{\noindent{\bf Proof of Theorem \ref{Thm. 4th powers}}.} In this proof, we set 
	$$D_4=\{b_1,\cdots,b_{(q-1)/4}\}=\{x^4:x\in\mathbb{F}_q^{\times}\}.$$
	It is easy to verify that 
	$$x^{(q-1)/4}-1=\prod_{1\le j\le (q-1)/4}\left(x-b_j\right).$$
	Noting that $q\equiv 5\pmod 8$, by the above we have 
	$$\prod_{1\le j\le (q-1)/4}b_j=(-1)^{(q+3)/4}=1.$$
	This implies 
	\begin{equation}\label{Eq. 1 in the proof of Thm. 2}
		\det(B_{q,4}(\chi^r))=\chi^r\left(\prod_{1\le i\le (q-1)/4}b_i\right)\cdot\det(N_q(\chi^r))=\det(N_q(\chi^r)),
	\end{equation}
	where 
	$$N_q(\chi^r)=\left[\chi^r\left(1+\frac{b_j}{b_i}\right)+\chi^r\left(1-\frac{b_j}{b_i}\right)\right]_{1\le i,j\le (q-1)/4}.$$
	
	For any integer $1\le k\le (q-1)/4$, let the column vector 
	$$\u_k:=\left(\chi^k(b_1),\cdots,\chi^k(b_{(q-1)/4})\right)^T,$$
	and let 
	\begin{equation}\label{Eq. 2 in the proof of Thm. 2}
		\mu_k=\sum_{1\le j\le (q-1)/4}\chi^r(1+b_j)\chi^k(b_j)+\sum_{1\le j\le (q-1)/4}\chi^r(1-b_j)\chi^k(b_j).
	\end{equation}
	Then for any $1\le i\le (q-1)/4$ one can verify that 
	\begin{align*}
		 &\sum_{1\le j\le (q-1)/4}\chi^r\left(1+\frac{b_j}{b_i}\right)\chi^k(b_j)+\sum_{1\le j\le (q-1)/4}\chi^r\left(1-\frac{b_j}{b_i}\right)\chi^k(b_j)\\
		=&\sum_{1\le j\le (q-1)/4}\chi^r\left(1+\frac{b_j}{b_i}\right)\chi^k(b_j/b_i)\chi^k(b_i)+\sum_{1\le j\le (q-1)/4}\chi^r\left(1-\frac{b_j}{b_i}\right)\chi^k(b_j/b_i)\chi^k(b_i)\\
		=&\sum_{1\le j\le (q-1)/4}\chi^r\left(1+b_j\right)\chi^k(b_j)\chi^k(b_i)+\sum_{1\le j\le (q-1)/4}\chi^r\left(1-b_j\right)\chi^k(b_j)\chi^k(b_i)\\
		=&\mu_k\cdot\chi^k(b_i).
	\end{align*}
	This implies that for any $1\le k\le (q-1)/4$ we have 
	$$N_q(\chi^r)\u_k=\mu_k\u_k.$$
	As $\u_1,\cdots,\u_{(q-1)/4}$ are linearly independent over $\mathbb{C}$, we see that 
	$\mu_1,\cdots,\mu_{(q-1)/4}$ are exactly all the eigenvalues of $N_q(\chi^r)$. We next evaluate the explicit value of $\mu_k$. 
	
	Let 
	$$h(x)=
	\frac{1}{4}\left(\varepsilon(x)+\chi^{\frac{q-1}{4}}(x)+\chi^{\frac{q-1}{2}}(x)+\chi^{\frac{3(q-1)}{4}}(x)\right).$$
	Then it is easy to see that 
	$$h(x)=\begin{cases}
		1 & \mbox{if}\ x\in D_4,\\
		0 & \mbox{otherwise.}
	\end{cases}$$
	By (\ref{Eq. 2 in the proof of Thm. 2}) one can verify that 
	\begin{equation}\label{Eq. 3 in the proof of Thm. 2}
	\mu_k
	=\sum_{x\in\mathbb{F}_q}h(x)\chi^r(1+x)\chi^k(x)+\sum_{x\in\mathbb{F}_q}h(x)\chi^r(1-x)\chi^k(x).
	\end{equation} 
	
	{\bf Case I:} $k\equiv 1\pmod 2$.
	
	In this case, by (\ref{Eq. 3 in the proof of Thm. 2}) for any $\psi_e\in\{\varepsilon,\chi^{\frac{q-1}{2}}\}$ (note that $\psi_e(-1)=1$) we have 
	\begin{align*}
		  &\sum_{x\in\mathbb{F}_q}\psi_e(x)\chi^r(1+x)\chi^k(x)+\sum_{x\in\mathbb{F}_q}\psi_e(x)\chi^r(1-x)\chi^k(x)\\
		 =&\sum_{x\in\mathbb{F}_q}\psi_e(x)\chi^r(1+x)\chi^k(x)+\sum_{x\in\mathbb{F}_q}\psi_e(-x)\chi^r(1+x)\chi^k(-x)\\
		 =&\sum_{x\in\mathbb{F}_q}\psi_e(x)\chi^r(1+x)\chi^k(x)-\sum_{x\in\mathbb{F}_q}\psi_e(x)\chi^r(1+x)\chi^k(x)\\
		 =&0.
	\end{align*}
	Also, for any $\psi_o\in\{\chi^{\frac{q-1}{4}},\chi^{\frac{3(q-1)}{4}}\}$ (note that $\psi_o(-1)=-1$) we have 
	\begin{align*}
		 &\sum_{x\in\mathbb{F}_q}\psi_o(x)\chi^r(1+x)\chi^k(x)+\sum_{x\in\mathbb{F}_q}\psi_o(x)\chi^r(1-x)\chi^k(x)\\
		=&\sum_{x\in\mathbb{F}_q}\psi_o(-x)\chi^r(1-x)\chi^k(-x)+\sum_{x\in\mathbb{F}_q}\psi_o(x)\chi^r(1-x)\chi^k(x)\\
		=&2\sum_{x\in\mathbb{F}_q}\psi_o(x)\chi^r(1-x)\chi^k(x)\\
		=&2J_q(\chi^r,\psi_o\chi^k).
	\end{align*}
    Hence 
    \begin{align*}
    	\lambda_k
    	=&\frac{1}{2}J_q(\chi^r,\chi^{\frac{q-1}{4}+k})+\frac{1}{2}J_q(\chi^r,\chi^{\frac{3(q-1)}{4}+k})\\
    	=&\frac{1}{2}J_q(\chi^r,\chi^{\frac{q-1}{4}+k})+\frac{1}{2}J_q(\chi^r,\phi\chi^{\frac{q-1}{4}+k})\\
    	=&\frac{1}{2}\sum_{x\in\mathbb{F}_q}\left(\varepsilon(x)+\phi(x)\right)\chi^{\frac{q-1}{4}+k}(x)\chi^r(1-x)\\
    	=&\frac{1}{2}\sum_{x\in\mathbb{F}_q}\chi^{\frac{q-1}{4}+k}(x^2)\chi^r(1-x^2)\\
    	=&\frac{1}{2}\sum_{x\in\mathbb{F}_q}\chi^{\frac{q-1}{2}+2k}(x)\chi^r(1-x)\chi^r(1+x)\\
    	=&(-1)^r\frac{q}{2}\cdot\ _{2}F_1\left(\begin{array}{cc}
    		\chi^{-r}  &  \chi^{2k+(q-1)/2}\\
    		&  \chi^{r+2k+(q-1)/2}
    	\end{array}\Bigg|-1\right)_q.
    \end{align*}
	The last equality follows from \cite[Definition 3.5]{Greene}.

	{\bf Case II:} $k\equiv 0\pmod 2$. 
	
	As in the case $k\equiv 1\pmod 2$, for any  $\psi_o\in\{\chi^{\frac{q-1}{4}},\chi^{\frac{3(q-1)}{4}}\}$ and $\psi_e\in\{\varepsilon,\chi^{\frac{q-1}{2}}\}$, by some calculations one can easily verify that 
	\begin{equation*}
		\sum_{x\in\mathbb{F}_q}\psi_o(x)\chi^r(1+x)\chi^k(x)+\sum_{x\in\mathbb{F}_q}\psi_o(x)\chi^r(1-x)\chi^k(x)=0,
	\end{equation*}
	and that 
	\begin{equation*}
		\sum_{x\in\mathbb{F}_q}\psi_e(x)\chi^r(1+x)\chi^k(x)+\sum_{x\in\mathbb{F}_q}\psi_e(x)\chi^r(1-x)\chi^k(x)=2J_q(\chi^r,\psi_e\chi^k).
	\end{equation*}
	Thus, we have 
	\begin{align*}
		 \lambda_k
		=&\frac{1}{2}J_q(\chi^r,\chi^k)+\frac{1}{2}J_q(\chi^r,\phi\chi^k)\\
		=&\frac{1}{2}\sum_{x\in\mathbb{F}_q}\chi^k(x^2)\chi^r(1-x^2)\\
		=&\frac{1}{2}\sum_{x\in\mathbb{F}_q}\chi^{2k}(x)\chi^r(1-x)\chi^r(1+x)\\
		=&(-1)^r\frac{q}{2}\cdot\ _{2}F_1\left(\begin{array}{cc}
			\chi^{-r}  &  \chi^{2k}\\
			&  \chi^{r+2k}
		\end{array}\Bigg|-1\right)_q.
	\end{align*}
	
	In view of the above, we obtain 
	\begin{align*}
		\det(B_{q,4}(\chi^r))
		=&\prod_{k=1}^{(q-5)/8}(-1)^r\frac{q}{2}\ _{2}F_1\left(\begin{array}{cc}
			\chi^{-r}  &  \chi^{4k}\\
			&  \chi^{r+4k}
		\end{array}\Bigg|-1\right)_q\\ \times
	    &\prod_{k=0}^{(q-5)/8}(-1)^r\frac{q}{2}\ _{2}F_1\left(\begin{array}{cc}
		\chi^{-r}  &  \chi^{4k+2+(q-1)/2}\\
		&  \chi^{r+4k+2+(q-1)/2}
	\end{array}\Bigg|-1\right)_q.
	\end{align*}
	This finally implies 
	\begin{equation*}
		\det(B_{q,4}(\chi^r))=(-1)^r\left(\frac{q}{2}\right)^{\frac{q-1}{4}}\cdot\prod_{k=0}^{(q-5)/4}\ _{2}F_1\left(\begin{array}{cc}
			\chi^{-r}  &  \chi^{4k}\\
			&  \chi^{r+4k}
		\end{array}\Bigg|-1\right)_q.
	\end{equation*}
	This completes the proof of Theorem \ref{Thm. 4th powers}.\qed

	\section{Proof of Theorem \ref*{Thm. Tq(r)}} 
		\setcounter{lemma}{0}
	\setcounter{theorem}{0}
	\setcounter{equation}{0}
	\setcounter{conjecture}{0}
	\setcounter{remark}{0}
	\setcounter{corollary}{0}

	In this section, we will prove our last theorem. Recall that $q=p^n$ is an odd prime power. Let $\zeta_{q-1}=e^{\frac{2\pi{\bf i}}{q-1}}$. Consider the cyclotomic field $K=\mathbb{Q}(\zeta_{q-1},\zeta_p)$. Let $\mathcal{O}_K$ be the ring of all algebraic integers over $K$ and let $\mathfrak{p}$ be a prime ideal of $\mathcal{O}_K$ with $p\in\mathfrak{p}$. Then it is known that 
	$$\mathcal{O}_K/\mathfrak{p}\cong\mathbb{F}_q.$$
	From now on, we identify $\mathbb{F}_q$ with $\mathcal{O}_K/\mathfrak{p}$. Let $\chi_{\mathfrak{p}}\in\widehat{\mathbb{F}_q^{\times}}$ be the Teich\"{u}muller character of $\mathfrak{p}$, i.e., 
	$$\chi_{\mathfrak{p}}(x)\mod\mathfrak{p}=x$$
	for any $x\in\mathbb{F}_q=\mathcal{O}_K/\mathfrak{p}$. It is easy to verify that $\chi_{\mathfrak{p}}$ is a generator of $\widehat{\mathbb{F}_q^{\times}}$. 
	
	We need the following known result (see \cite[Proposition 3.6.4]{Cohen}).
	
	\begin{lemma}\label{Lem. congruence for Jacobi sums}
		Let $1\le a,b\le q-2$ be integers. Then 
		$$J_q(\chi_{\mathfrak{p}}^{-a},\chi_{\mathfrak{p}}^{-b})\equiv -\binom{a+b}{a}=-\frac{(a+b)!}{a!b!}\pmod{\mathfrak{p}},$$
		and 
			$$J_q(\chi_{\mathfrak{p}}^{-a},\chi_{\mathfrak{p}}^{0})\equiv -\binom{a+0}{a}\equiv -1\pmod{\mathfrak{p}}.$$
		In particular, if $a+b\ge q$, then 
		$$J_q(\chi_{\mathfrak{p}}^{-a},\chi_{\mathfrak{p}}^{-b})\equiv 0\pmod{\mathfrak{p}}.$$
	\end{lemma}
	
	Now we are in a position to prove our last theorem.
	
	{\noindent{\bf Proof of Theorem \ref{Thm. Tq(r)}}.} (i) Suppose $(q-1)/k\equiv 0\pmod 2$. Then by the definition of the Teich\"{u}muller character $\chi_{\mathfrak{p}}$ and Theorem \ref{Thm. Bq(r)}(ii) we have   
	\begin{equation*}
		\det(T_{q,k}(r))=\det(B_{q,k}(\chi^{-r}))\mod{\mathfrak{p}}=0\mod{\mathfrak{p}}.
	\end{equation*}
	This completes the proof of (i).
	
	(ii) As  
	\begin{equation*}
		\det(T_{q,2}(r))=\det(B_{q,2}(\chi^{-r}))\mod{\mathfrak{p}},
	\end{equation*}
    by Theorem \ref{Thm. Bq(r)}(ii) we see that $T_{q,2}(r)$ is singular whenever $q\equiv 1\pmod 4$. Below we suppose $q\equiv 3\pmod 4$. By Theorem \ref{Thm. Bq(r)} and Lemma \ref{Lem. congruence for Jacobi sums} we obtain 
	\begin{align}\label{Eq. det Tq(r) as the product of binomial coeffients}
		  \det(T_{q,2}(r))
		&=\prod_{0\le k\le (q-3)/2}J_q(\chi_{\mathfrak{p}}^{-r},\chi_{\mathfrak{p}}^{-2k}) \mod{\mathfrak{p}} \notag\\
		&=(-1)^{\frac{q-1}{2}}\prod_{0\le k\le (q-3)/2}\binom{r+2k}{r} \mod{\mathfrak{p}}.
	\end{align}
	By the above, it is easy to see that 
	$$T_{q,2}(r)\ \text{is nonsingular}\ \Leftrightarrow \binom{r+2k}{r}\not\equiv 0\pmod p\ \text{for any}\ 0\le k\le (q-3)/2.$$
	Suppose $q-2\ge r\ge3$. Then by Lemma \ref{Lem. congruence for Jacobi sums} we have 
	$$\binom{r+q-3}{r}\equiv 0\pmod p$$ and hence $T_{q,2}(r)$ is singular. Suppose now $r\in\{1,2\}$. If $n\ge 2$, then $(p-1)/2\le (q-3)/2$ and hence 
	$$\binom{r+p-1}{r}\equiv 0\pmod p.$$
	This implies that $T_{q,2}(r)$ is singular. Suppose now $r\in\{1,2\}$ and $n=1$. Then for any $0\le k\le (p-3)/2$ we clearly have
	$$\binom{r+2k}{r}\not\equiv 0\pmod p.$$
	Hence $T_{q,2}(r)$ is nonsingular. 
	
	Suppose now $p\equiv 3\pmod 4$ and $p>3$. Then by (\ref{Eq. det Tq(r) as the product of binomial coeffients}) one can verify that 
	\begin{align*}
		 \det(T_{p,2}(1))
		&\equiv -\prod_{0\le k\le (p-3)/2}(2k+1) \\
		&\equiv \frac{-(p-1)!}{2^{(p-1)/2}\cdot\frac{p-1}{2}!} \\
		&\equiv (-1)^{\frac{p+1}{4}+\frac{h(-p)+1}{2}} \pmod{\mathfrak{p}}.
	\end{align*}
	The last equality follows from $2^{(p-1)/2}\equiv (-1)^{(p+1)/4} \pmod p$ and the Mordell congruence \cite{Mordell}
	$$\frac{p-1}{2}!\equiv (-1)^{\frac{h(-p)+1}{2}} \pmod p.$$
	 For $r=2$, by (\ref{Eq. det Tq(r) as the product of binomial coeffients}) we have 
	\begin{align*}
		\det (T_{p,2}(2))
		&\equiv -\prod_{0\le k\le (p-3)/2}(k+1)(2k+1)\\
		&\equiv -\frac{p-1}{2}!\prod_{0\le k\le (p-3)/2}(2k+1)\\
		&\equiv -\frac{p-1}{2}!\cdot\frac{(p-1)!}{2^{(p-1)/2}\cdot\frac{p-1}{2}!}\\
		&\equiv (-1)^{\frac{p+1}{4}} \pmod{\mathfrak{p}}.
	\end{align*}
	The last equality follows from $2^{(p-1)/2}\equiv (-1)^{(p+1)/4} \pmod p$. This completes the proof of (ii).
	
	(iii) Suppose $q\equiv 5\pmod 8$. Then by Theorem \ref{Thm. 4th powers}, Remark \ref{Remark of Thm. 2} and Lemma \ref{Lem. congruence for Jacobi sums} we have 
	\begin{align*}
	\det (T_{q,4}(r))
	&=\det(B_{q,4}(\chi_{\mathfrak{p}}^{-r}))\mod{\mathfrak{p}}\\
	&=\left(\frac{1}{2}\right)^{\frac{q-1}{4}}\prod_{0\le k\le (q-5)/4}\left(J_q(\chi_{\mathfrak{p}}^{-r},\chi_{\mathfrak{p}}^{-2k})+J_q(\chi_{\mathfrak{p}}^{-r},\chi_{\mathfrak{p}}^{-2k-(q-1)/2})\right) \mod{\mathfrak{p}}\\
	&=\left(\frac{1}{2}\right)^{\frac{q-1}{4}}\prod_{k=0}^{(q-5)/4}
	\left(\binom{2k+r}{r}+\binom{2k+r+(q-1)/2}{r}\right)\mod{\mathfrak{p}}\in\mathbb{F}_p.
	\end{align*}

	In view of the above, we have completed the proof of Theorem \ref{Thm. Tq(r)}.\qed 
	
	\Ack\ The authors would like to thank the referee for helpful comments.
	
	This work was supported by the Natural Science Foundation of China (Grant Nos. 12101321 and 12201291). The first author was also supported by Natural Science Foundation of Nanjing University of Posts and Telecommunications (Grant No. NY224107).

\end{document}